\documentclass[12pt]{article}

\usepackage{amsfonts,graphicx}
\usepackage{amsmath}
\usepackage{amssymb}
\usepackage{rotating}
\usepackage{enumerate} 
\usepackage{framed}
\usepackage{amsrefs}
\usepackage{amsthm}

\renewcommand{\labelenumi}{\roman{enumi})}
\newtheorem{thm}{Theorem}[section]
\newtheorem{prop}[thm]{Proposition}
\newtheorem{lem}[thm]{Lemma}
\newtheorem{cor}[thm]{Corollary}
\numberwithin{equation}{section}

\newcommand{\n}{\Vert}

\newcommand{\IIi}{${\textrm{II}}_1\ $}
\newcommand{\ce}[2]{\mathbb E_{#1}\left(#2\right)}
\newcommand{\nrm}{\mathcal{N}_M}
\newcommand{\reli}{N'\cap \langle M, e_N \rangle}
\newcommand{\tr}{{\textrm{Tr}}}
\newcommand{\mi}{\langle M,e_N\rangle}

\newcommand{\vpic}[2]{\mbox{$\begin{array}[c]{l}\includegraphics[width=#2]{#1}
\end{array}$}}
\title{Strong Singularity for Subfactors}

\author{Pinhas Grossman\thanks{Partially supported by NSF Grant DMS-0801235}\\\normalsize\texttt{pinhas.grossman@vanderbilt.edu}
\and Alan Wiggins\\ 
\normalsize\texttt{alan.d.wiggins@vanderbilt.edu}}

\date{}

\begin{document}
\maketitle
\begin{abstract}
We examine the notion of $\alpha$-strong singularity for subfactors of a \IIi factor, which is a metric quantity that relates the distance between a unitary in the factor and a subalgebra with the distance between that subalgebra and its unitary conjugate. Through planar algebra techniques, we demonstrate the existence of a finite index singular subfactor of the hyperfinite \IIi factor that cannot be strongly singular with $\alpha=1$, in contrast to the case for masas. Using work of Popa, Sinclair, and Smith, we show that there exists an absolute constant $0<c<1$ such that all singular subfactors are $c$-strongly singular. Under the hypothesis of $2$-transitivity, we prove that finite index subfactors are $\alpha$-strongly singular with a constant that tends to $1$ as the Jones Index tends to infinity and infinite index subfactors are $1$-strongly singular. Finally, we give a proof that proper finite index singular subfactors do not have the weak asymptotic homomorphism property relative to the containing factor.
  
\end{abstract}

\section{Introduction}

The study of subfactors of a \IIi factor was initiated by Vaughan Jones in \cite{Jones.Index}, where he defined the index $[M:N]$ of a subfactor inclusion $N\subseteq M$ to be dimension of $L^2(M)$ as a left Hilbert $N$-module. He showed that the spectrum of possible index values contains both a continuous and a discrete part: while the index can take any value greater than or equal to $4$, values less than $4$ are necessarily of the form $4\cos^2(\frac{\pi}n)$ for an integer $n \geq 3$. He also showed that all admissible index values are realized by subfactors of the hyperfinite \IIi factor.

Since the introduction of the index, invariants of subfactors both numerical and otherwise have been a rich area of study. It was shown in \cite{Jones.Index} that if $[M:N]<\infty$, then $N'\cap \langle M, e_N\rangle$ is finite dimensional, and by repeating the Jones construction one obtains a double sequence of inclusions of finite dimensional algebras 
\[
\begin{array}{ccccccc} N'\cap M & \subseteq & \reli & \subseteq & N'\cap \langle \langle M,e_N\rangle , e_M \rangle & \subseteq & \dots \\ \begin{sideways}$\subseteq$\end{sideways} & & \begin{sideways}$\subseteq$\end{sideways} & & \begin{sideways}$\subseteq$\end{sideways} & & \\\mathbb{C}I & \subseteq & M'\cap \mi & \subseteq & M'\cap \langle \langle M,e_N\rangle , e_M \rangle & \subseteq & \dots \end{array}
\]
called the standard invariant. Popa showed in \cite{Popa.CBMSNotes} that when $[M:N]\leq 4$ and $M$ is hyperfinite, the subfactor can be reconstructed from the standard invariant. He also proved that this holds in general for a larger class of \emph{strongly amenable} subfactors. In \cite{Dietmar.Index6} it was shown that there exist infinitely many nonisomorphic inclusions of index 6 subfactors of the hyperfinite \IIi factor with the same standard invariant. It is a major open question in subfactor theory to decide whether all possible standard invariants occur for subfactors of the hyperfinite \IIi factor. A subproblem, still open, is whether all index values can be obtained for irreducible subfactors, those with $N'\cap M=\mathbb{C}I$. Singular subfactors fall under the umbrella of this latter problem. 

The concept of singularity for a subalgebra of a \IIi factor $M$ dates back to Jacques Dixmier \cite{Dixmier.Masa} in the context of maximal abelian *-subalgebras (masas) $A$ of $M$. If $\mathcal{U}(M)$ is the group of unitaries in $M$, and
\[
\nrm (A)=\{u\in \mathcal{U}(M):uAu^*=A\},
\]
$A$ is said to be singular if $\nrm(A)=\mathcal{U}(A)$. Dixmier provided examples of singular masas in the hyperfinite \IIi factor $R$. 

In general, it is difficult to tell whether a given masa is singular or not. This difficulty led Sinclair and Smith to define the notion of $\alpha$-strong singularity in \cite{Sinclair.strongsing} as an analytic quantity which would imply the algebraic condition of singularity for masas in a \IIi factor $M$. The definition was extended to subalgebras $B$ in \cite{Sinclair.StrongSing2} by Sinclair, Smith, and Robertson. 

In \cite{singular.al}, it was shown that singularity and strong singularity for masas are equivalent to a formally stronger property which first appeared in \cite{Sinclair.StrongSing2} and has come to be known as the weak asymptotic homomorphism property, or WAHP, in $M$. Using the equivalence between the WAHP and singularity, it was shown in \cite{singular.al} that the tensor product of singular masas in \IIi factors is again a singular masa in the tensor product factor. 

It is natural to ask what relationships these properties have for arbitrary subalgebras of $M$. Herein, we will consider the case where $B=N$ is a subfactor of $M$. Understanding singular subfactors of the hyperfinite \IIi factor would shed light on (and possibly decide) whether all index values occur for irreducible subfactors for subfactors of the hyperfinite \IIi factor. It is worth noting that all ``small'' noninteger index values (between 3 and 4) yield singular subfactors.

One might guess that in this situation, the other extreme from masas, that it would be possible to deduce strong singularity from singularity. The main result of this paper is that this is not the case. Using planar algebra techniques, we produce an example of a finite index singular subfactor of the hyperfinite \IIi factor $\mathcal{R}$. that is no more than $\sqrt{2({\sqrt{2}-1})}$-strongly singular in $\mathcal{R}$. The supremum over all admissible numbers $\alpha$ appearing in the strong singularity inequality then represents a new, nontrivial numerical invariant for singular subfactors of \IIi factors under unitary conjugacy. 

The paper is organized as follows: Section \ref{prelims} establishes notation and general background. In Section \ref{mainstuff}, we establish positive results for strong singularity constants. Theorem \ref{relativecomm} shows that when the higher relative commutant $N' \cap \langle M,e_N\rangle$ is $2$-dimensional, proper finite index subfactors of $M$ are $\alpha$-strongly singular in $M$ where $\alpha=\sqrt{\frac {[M:N]-2}{[M:N]-1}}$.  As this constant tends to one as the index tends to infinity, this suggests that infinite index singular subfactors are strongly singular. Indeed, the methods of proof for the finite index case yield strong singularity for an infinite index inclusion $N\subseteq M$ when $\reli$ is $2$-dimensional. Using results from \cite{Sinclair.PertSubalg}, we obtain an absolute constant $c=\frac 1 {13}$ for which all singular subfactors are $c$-strongly singular. 

In Section \ref{counterex}, we give the example described above as the unique (up to conjugacy) subfactor of index $2+\sqrt{2}$ subfactor of the hyperfinite \IIi factor, and prove the aforementioned upper bound on $\alpha$. Using Theorem \ref{relativecomm}, we can obtain a lower bound of $\sqrt{2-\sqrt{2}}$. Finally, we give a simple proof in  Section \ref{nowahp} of the fact, due to Popa, that no proper singular finite index subfactor of $M$ has the WAHP. Thus, these properties cannot be equivalent in general. 

Let us briefly discuss existence questions for singular subfactors. Since technically $M$ is a strongly singular subfactor of itself (with the WAHP), existence questions for singular (or strongly singular) subfactors must be qualified. Recently, Stefaan Vaes has proved that there exists a factor $M$ such that every finite index irreducible subfactor is equal to $M$ \cite{Vaes.nontrivialsubfactor}. This implies that there exist factors with no proper finite index singular subfactors. On the other hand, Popa has shown in \cite{Popa.SingularMasas} that there always exist singular masas in separable \IIi factors. The correct analog of the question for masas, then, is to ask whether there always exist infinite index hyperfinite singular or $\alpha$-strongly singular subfactors of any separable \IIi factor. 

An example in the hyperfinite \IIi factor of an infinite index subfactor with the WAHP was provided in \cite{Sinclair.StrongSing2}. In \cite{Popa.MaxInjective}, Popa remarks that by results from \cite{Popa.Kadison}, every separable \IIi factor has a semi-regular masa that is contained in some (necessarily irreducible) hyperfinite subfactor, and so by Zorn's Lemma has an irreducible maximal hyperfinite subfactor. Such an object is then a maximal hyperfinite subalgebra of $M$, and as Popa observes, any maximal hyperfinite subalgebra is singular \cite{Popa.MaxInjective}. Therefore, any separable \IIi factor has an infinite index hyperfinite singular subfactor. Maximal hyperfinite subfactors in any \IIi factor were first exhibited in \cite{Fuglede.MaxInjective}. Whether there exist strongly singular hyperfinite subfactors or hyperfinite subfactors with the WAHP in any \IIi factor remains an open question. 

\section{Preliminaries and Notation}\label{prelims}

Throughout, $M$ will denote a \IIi factor and $N$ a subfactor of $M$. Unless otherwise noted, $M$ shall be regarded as faithfully represented on the Hilbert space $L^2(M)=L^2(M, \tau)$, where $\tau$ denotes the unique normal, faithful, tracial state on $M$. Elements of $M$ considered as a subspace of $L^2(M)$ shall be denoted by $\hat{x}$ or $x\hat{I}$ for $x$ in $M$ and $I$ the identity element of $M$. The element $\hat{I}$ is a cyclic and separating vector for $M\subseteq B(L^2(M))$. If $J$ denotes the isometric involution on $L^2(M)$ defined by 
\begin{equation}
J(x\hat{I})= x^*\hat{I}
\end{equation}
on $M\hat{I}$, then $JMJ=M'$. 

If $e_N$ is the orthogonal projection from $L^2(M)$ onto $L^2(N)$, then the von Neumann algebra $\langle M, e_N\rangle$ generated by $M$ and $e_N$ is a factor of type II equal to $JN'J$, and so is of type \IIi if and only if $N'$ is. We will denote by $\mathbb{E}_N$ the unique normal, faithful, trace-preserving conditional expectation from $M$ onto $N$, which can be thought of as the restriction of $e_N$ to $M\hat{I}$. The factor $\langle M, e_N\rangle$ possesses a unique normal, faithful, semifinite tracial weight $\tr$ such that for all $x,y$ in $M$,
\begin{enumerate}
\item $\tr(I)=[M:N]$;
\item $e_N\langle M,e_N \rangle e_N=Ne_N$;
\item $\tr (xe_Ny)=\tau(xy)$;
\item $e_Nxe_N=\ce{N}{x}e_N=e_N\ce{N}{x}$.
\end{enumerate}
If $[M:N]<\infty$, then for every element $x$ in $\langle M, e_N \rangle$, there is a unique element $y$ in $M$ with $xe_N=ye_N$. Proofs of these facts may be found in \cite{Jones.SubfactorsBook} or \cite{Jones.Index}. We shall denote by Aut$(N)$, $\mathcal{U}(N)$, and $\nrm(N)$ the groups of automorphisms of $N$, unitaries in $N$, and normalizing unitaries of $N$ in $M$, respectively.

A von Neumann subalgebra $B$ of $M$ is $\alpha$-strongly singular \cite{Sinclair.strongsing} if there is a constant $0< \alpha\leq 1$ such that for all unitaries $u\in M$, 
\begin{equation}\label{ss}
\alpha \n u-\ce{B}{u}\n_2 \leq \n \mathbb{E}_B-\mathbb{E}_{uBu^*}\n_{\infty, 2}
\end{equation}
where $\n T\n _{\infty,2}=\displaystyle\sup_{\genfrac{}{}{0cm}{1}{x\in M}{\n x\n \leq 1}}\n Tx\n_2.$ If $\alpha=1$, then $B$ is said to be strongly singular. 

\section{Positive Strong Singularity Results for Subfactors}\label{mainstuff}

Under the assumption that $N'\cap \langle M,e_N\rangle $ is 2-dimensional, we can establish strong singularity constants for singular subfactors. In \cite{G.J}, this condition is referred to as 2-transitivity. Note that if $[M:N]>2$ and $N'\cap \langle M,e_N\rangle$ is 2-dimensional, then $N$ is automatically singular in $M$, as any $u$ in $\mathcal{N}_M(N)\backslash \mathcal{U}(N)$ yields the projection $ue_Nu^*$ in $N'\cap \langle M, e_N\rangle$. This projection is not $e_N$ since $\{e_N\}'\cap M=N$ and it is also not $e_N^{\perp}$ since 
\[
\tr (e_N^{\perp})>\tr (e_N)=\tr(ue_Nu^*).
\]
By Goldman's Theorem (\cite{Goldman.Index2} or \cite{Jones.SubfactorsBook}), all index 2 subfactors are regular. Any subfactor of index strictly between 3 and 4 is 2-transitive, and therefore singular. There also exist 2-transitive hyperfinite subfactors for every integer index $\geq 3 $. 

\renewcommand{\labelenumi}{\arabic{enumi})}
\begin{thm}\label{relativecomm}
Let $N\subseteq M$ be a singular subfactor with $N'\cap \langle M,e_N\rangle$ 2-transitive. If $[M:N]<\infty$, then $N$ is $\sqrt{\frac {[M:N]-2}{[M:N]-1}}$-strongly singular in $M$. If $[M:N]=\infty$, then $N$ is strongly singular in $M$.
\end{thm}
\begin{proof} Let $N\subseteq M$ be a singular inclusion of subfactors and suppose $N'\cap \langle M, e_N \rangle $ is 2-dimensional. Let $u$ be a unitary in $M$ and define $C$ to be the weakly closed convex hull of the set $\{we_Nw^* : w\in uNu^*\}$ where $w$ is unitary. Then $C$ admits a unique element of minimal $2$-norm denoted by $h$ which has the following properties, detailed in \cite{Sinclair.PertSubalg} :
\begin{enumerate}
\item $h\in (uNu^*)'\cap \langle M, e_N\rangle$;

\item Tr$(e_Nh)={\textrm{Tr}}(h^2)$;

\item Tr$(h)=1$;

\item $1-{\textrm{Tr}} (e_Nh)\leq \n \mathbb{E}_N-\mathbb{E} _{uNu^*}\n ^2 _{\infty, 2}$.

\end{enumerate}

Now $(uNu^*)'\cap \langle M, e_N \rangle$ has basis $ue_Nu^*$ and $ue_N^{\perp}u^*$, so that $h=\alpha (ue_Nu^*)+\beta (ue_N^{\perp}u^*)$ for some scalars $\alpha$ and $\beta$. By 3), 
\[
1=\tr (h)=\alpha +\lambda \beta,
\]
where $\lambda=[M:N]-1$. If $[M:N]=\infty$, then $\tr (e_N ^{\perp})=\infty$, which implies that $\beta =0$, and therefore $\alpha=1$. If $[M:N]$ is finite, then $\alpha=1-\lambda \beta$, and expanding $\tr (e_Nh)$ yields
\begin{align*}
{\textrm{Tr}}(e_Nh)=&\alpha \tr (e_Nue_Nu^*)+\beta \tr(e_N-e_Nue_Nu^*)\\
=&\alpha \tr (e_N\mathbb{E}_N(u)\mathbb{E}_N(u^*))+\beta (\tr(e_N)-\tr (e_N\mathbb{E}_N (u)\mathbb{E}_N (u^*)))\\
=&\alpha \tau (\mathbb{E}_N (u)\mathbb{E}_N (u^*))+\beta (1-\tau (\mathbb{E}_N(u) \mathbb{E}_N (u^*)))\\
=&\alpha \n \mathbb{E}_N(u)\n^2 _2+\beta \n u-\mathbb{E}_N(u)\n ^2 _2,
\end{align*}
since $1=\n u\n^2 _2= \n \mathbb{E}_N(u)\n^2 _2+\n u-\mathbb{E}_N(u)\n^2 _2$. Setting $k=\n u-\mathbb{E}_N(u)\n^2 _2$ and substituting the formula for $\alpha$ gives 
\[
\tr (e_Nh)=(1-\lambda \beta)(1-k)+\beta k.
\]
Using 2), we have
\[
(1-\lambda \beta)(1-k)+\beta k=\tr (e_Nh)=\tr (h^2)=(1-\lambda \beta)^2+\lambda \beta ^2,
\]
and so $(\lambda^2+\lambda)\beta^2-(\lambda +k+\lambda k)\beta +k=0$. We may then solve for $\beta$ in terms of $\lambda$ and $k$, obtaining the roots $\beta=\dfrac k \lambda$ and $\beta =\dfrac 1 {1+\lambda}$.

Suppose that $\beta=\dfrac 1 {1+\lambda}$. Then 
\[
\alpha=1-\frac {\lambda}{1+\lambda}=\frac 1 {1+\lambda}=\beta,
\]
 and so 
\[
h=\beta (ue_Nu^*)+\beta (ue_N^{\perp}u^*)=\beta I=\dfrac 1 {1+\lambda}I.
\]
 Since $h$ is an element of $C$, there exist natural numbers $\{n_j\}_{j=1} ^{\infty}$, positive reals $\{\gamma_i ^{(j)}\}_{i=1} ^{n_j}$ with $\sum_{i=1}^{n_j}\gamma_i^{(j)}=1$ and unitaries $\{w_i ^{(j)}\}_{i=1} ^{n_j}$ in $N$ with 
\[
\lim_{j\to \infty}\sum_{i=1} ^{n_j}\gamma_i^{(j)}uw_i^{(j)}u^*e_Nu(w_i^{(j)})^*u^*=\frac 1 {1+\lambda} I
\]
 in WOT. 

Then also 
\[
\lim_{j\to \infty}\sum_{i=1} ^{n_j}\gamma_i^{(j)}w_i^{(j)}u^*e_Nu(w_i^{(j)})^*= \frac 1 {1+\lambda} I
\]
 in WOT and 
 \[
 \lim_{j\to \infty}e_N\left(\sum_{i=1} ^{n_j}\gamma_i^{(j)}w_i^{(j)}u^*e_Nu(w_i^{(j)})^*\right)= \frac {e_N} {1+\lambda}
 \] 
 in WOT. Taking the trace of both sides yields 
 \[
 \tr \left( \lim_{j\to \infty}e_N\left(\sum_{i=1} ^{n_j}\gamma_i^{(j)}w_i^{(j)}u^*e_Nu(w_i^{(j)})^*\right)\right)= \tr \Big{(}\frac {e_N} {1+\lambda}\Big{)}
=\frac 1 {1+\lambda}.
\]
 However, for any $n_j$, $1\leq j< \infty$,
\begin{align*}
\tr \left(e_N\left(\sum_{i=1} ^{n_j}\gamma_i^{(j)}w_i^{(j)}u^*e_Nu(w_i^{(j)})^*\right)\right)&= \tr \left(e_N\left(\sum_{i=1} ^{n_j}\gamma_i^{(j)}w_i^{(j)}\ce{N}{u^*}\ce{N}{u}(w_i^{(j)})^*\right)\right)\\
&=\tau \left(\sum_{i=1} ^{n_j}\gamma_i^jw_i^{(j)}\ce{N}{u^*}\ce{N}{u}(w_i^{(j)})^*\right)=\n \ce{N}{u}\n^2 _2,
\end{align*}
and so $\dfrac 1{1+\lambda} =\n \ce{N}{u}\n^2 _2$. We obtain that 
\[
k=1-\n \ce{N}{u}\n^2 _2=1-\dfrac 1{1+\lambda}=\dfrac{\lambda}{1+\lambda}.
\]
Then $\dfrac k \lambda=\dfrac 1 {1+\lambda}$, and so the only instance where $\beta=\dfrac 1 {1+\lambda}$ is when $k=\dfrac {\lambda}{1+\lambda}$, and there the two roots are identical.

We may then take $\beta=\dfrac k \lambda$ and so $\alpha=1-\lambda\beta=1-k$ when $[M:N]$ is finite. Hence
\begin{equation}
h=(1-k)(ue_Nu^*)+\dfrac k \lambda (ue_N^{\perp}u^*). 
\end{equation}
By 4),
\[
\n \mathbb{E}_N-\mathbb{E} _{uNu^*}\n ^2 _{\infty, 2}\geq 1-\tr(e_Nh)=1-\Big{(}(1-k)^2+\frac{k^2}={\lambda}\Big{)}=k\Big{(}2-\Big{(}1+\frac 1 \lambda \Big{)}k\Big{)}
\]
and therefore 
\begin{equation}\label{quad}
\n u-\ce{N}{u}\n^2_2=k\leq \dfrac 1{2-\Big{(}1+\frac 1 \lambda\Big{)}k}\n \mathbb{E}_N-\mathbb{E} _{uNu^*}\n ^2 _{\infty, 2}.
\end{equation} 
As $k\leq 1$, 
\[
2-\Big(1+\dfrac 1 \lambda \Big{)}k\geq 2-\Big{(}1+\dfrac 1 \lambda \Big{)}=1-\dfrac 1 \lambda,
\]
and it follows that
\begin{align*}
\n u-\ce{N}{u}\n^2_2\leq \frac 1 {1-\frac 1 \lambda}\n \mathbb{E}_N-\mathbb{E} _{uNu^*}\n ^2 _{\infty, 2}&=\frac \lambda {\lambda -1}\n \mathbb{E}_N-\mathbb{E} _{uNu^*}\n ^2 _{\infty, 2}\\
&=\frac {[M:N]-1}{[M:N]-2}\n \mathbb{E}_N-\mathbb{E} _{uNu^*}\n ^2 _{\infty, 2}.
\end{align*}
Hence $N$ is $\displaystyle \sqrt{\frac {[M:N]-2}{[M:N]-1}}$-strongly singular in $M$.

If $[M:N]=\infty$, then as previously noted, $\alpha=1$ and so $h=ue_Nu^*$. Therefore, 
\[
\n u-\ce{N}{u}\n^2_2=1-\tr (e_Nh)\leq \n \mathbb{E}_N-\mathbb{E} _{uNu^*}\n ^2 _{\infty, 2}
\]
so that $N$ is strongly singular in $M$ and the proof is complete. \end{proof}

In the situation of Theorem \ref{relativecomm}, we may immediately show that when unitaries are close to a finite index singular subfactor in $2$-norm, they satisfy the equation for strong singularity.

\begin{cor} Under the hypotheses of Theorem \ref{relativecomm}, if $[M:N]<\infty$ and 
\[
\n u-\ce{N}{u}\n_2\leq \sqrt{\frac {[M:N]-1} {[M:N]}},
\]
 then 
\[
\n u-\ce{N}{u}\n_2\leq \n\mathbb{E}_N-\mathbb{E}_{uNu^*}\n_{\infty,2}.
\]
\end{cor}
\begin{proof} Recall $k=\n u-\ce{N}{u}\n_2^2$ and $\lambda=[M:N]-1$. If $k\leq \frac {\lambda} {\lambda +1}$, then 
\[
2-\left(1+\frac 1 \lambda\right)k=2-\left(\frac {\lambda +1}{\lambda}\right)k\geq 2-\left(\frac {\lambda +1}{\lambda}\right)\left(\frac {\lambda}{\lambda+1}\right)=1
\]
Using equation \eqref{quad},
\[
\n u-\ce{N}{u}\n_2^2=k\leq k\left(2-\left(1+\frac 1 \lambda\right)k\right)\leq \n \mathbb{E}_N-\mathbb{E} _{uNu^*}\n ^2 _{\infty, 2}.
\]
\end{proof}

We end this section by producing an absolute constant $\alpha$ for which all singular subfactors of $M$ are strongly singular. First, we need a lemma dealing with the form of certain partial isometries. This fact is well-known, but we include a proof for completeness.

\begin{lem}\label{partialiso} (Dye's Theorem for subfactors) Let $M$ be a {\rm \IIi}factor and let $N$ be a subfactor of $M$. Suppose $v$ is a partial isometry in $M$ with $vNv^*\subseteq N$ and $v^*v\in N$. Then there exists a unitary $u\in M$ with $uNu^*\subseteq N$ and $v=uv^*v$. If $v^*Nv\subseteq N$ as well, then $u$ can be chosen to be in $\nrm (N)$.
\end{lem}
\begin{proof}

Suppose $v^*v\ne 0$. Let $n$ be the integer such that $(n-1)\cdot\tau(p)<1\leq n\cdot\tau(p)$. Let $p_1=v^*v$ and let $\{p_i\}_{i=2}^{n-1}\subseteq N$ be pairwise orthogonal projections subordinate to $p_1^{\perp}$ with $\tau(p_i)=\tau(p_1)$. Set $p_n=I-\sum_{i=1}^{n-1}p_i$. 

For each $2\leq i\leq n$ take $w_i$ to be a partial isometry in $N$ with $w_i^*w_i=p_i$ and $w_iw_i^*\leq p_1$. Set $w_1=p$. Then the projections $vw_ip_iw_i^*v^*\in N$ for all $1\leq i\leq n$. With $q_1=vv^*$, choose a collection of pairwise orthogonal projections $\{q_i\}_{i=2}^{n}\subseteq N$ in the same manner as for $p_1$. Take partial isometries $\{v_i\}_{i=1}^n\subseteq N$ so that $v_1=q$, $v_iv_i^*=q_i$ and $v_i^*v_i=vw_iw_i^*v^*\leq q_1$ for $2\leq i\leq n$. Then $\displaystyle u=\sum_{i=1}^nv_ivw_i$ is a unitary in $M$, and if $x\in N$, 
\[
uxu^*=\left(\sum_{i=1}^nv_ivw_i\right)x\left(\sum_{j=1}^nw_j^*v^*v_j^*\right)=\sum_{i,j=1}^nv_iv(w_ixw_j^*)v^*v_j\in N
\]
 as $w_ixw_j^*\in N$, $vNv^*\subseteq N$, and $v_i\in N$ for all $1\leq i,j\leq n$. Therefore $uNu^*\subseteq N$. Since $w_ip_1=\delta_{i,1}p_1$, we obtain that $v=up_1=uv^*v$.

We have shown that if $vNv^*\subseteq N$ and $v^*v\in N$, then $v=uv^*v$ for some unitary $u\in M$ with $uNu^*\subseteq N$. If in addition $v^*Nv\subseteq N$, then one can check that $u^*Nu\subseteq N$, so that $u$ normalizes $N$. \end{proof}

The following result is Theorem 5.4 in \cite{Sinclair.PertSubalg}.
\begin{thm}\label{pert}{\rm(Popa, Sinclair, \& Smith)}
Suppose $\delta >0$ and $N$, $N_0$ are two subfactors of $M$ with $\n \mathbb{E}_{N}-\mathbb{E}_{N_0}\n_{\infty, 2}\leq\delta$. Then there exist projections $q_0\in N_0$, $q\in N$, $q_0'\in N_0'\cap M$, $q'\in N'\cap M$, $p_0=q_0q_0'$, $p=qq'$, and a partial isometry $v$ in $M$ such that $vp_0N_0p_0v^*=pNp$, $vv^*=p$, $v^*v=p_0$, and 
\begin{equation}\label{isom}
\n 1-v\n_2\leq 13\delta, \quad \tau(p)=\tau(p_0)\geq 1-67\delta^2.
\end{equation}
\end{thm}
There is a similar theorem for arbitrary subalgebras of $M$, also in \cite{Sinclair.PertSubalg}.  As a direct consequence of Theorem \ref{pert} and Lemma \ref{partialiso}, we obtain

\begin{cor} Let $N$ be a singular subfactor in $M$. Then $N$ is $\frac 1 {13}$-strongly singular in $M$. \end{cor}

Theorem \ref{relativecomm} shows that the constant $\frac 1 {13}$ is not always optimal even for proper finite index singular subfactors.

\section{A Singular Subfactor that is Not Strongly Singular}\label{counterex}

In this section we describe an example, suggested to the authors by Vaughan Jones, of a subalgebra of the hyperfinite II$_1$ factor which is singular but not strongly singular. The subalgebra is in fact the unique (up to conjugacy) subfactor with index $2+\sqrt{2}$. To show it is not strongly singular, we will estimate both sides of inequality \ref{ss} for a specific unitary.

In \cite{G.J}  an irreducible quadrilateral of hyperfinite II$_1$ factors 
$\begin{array}{ccc}
P & \subset& M \\
\cup & & \cup \\
N &\subset &Q 
\end{array}$ was constructed such that each of the four elementary subfactors $$N \subset P, N \subset Q, P \subset M, Q \subset M $$ has index $2+\sqrt{2}$. The principal graph of $N\subseteq M$ is given as\\
\includegraphics[width=2in]{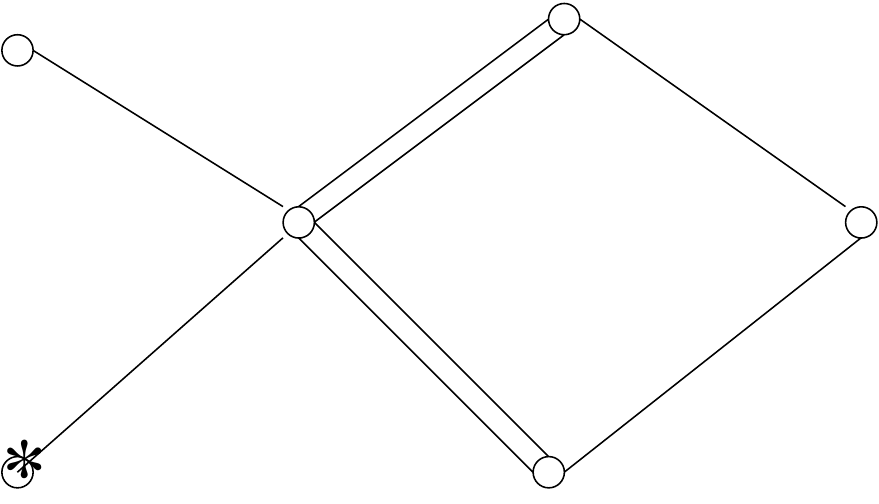}
\\

It is shown in \cite{G.J} that such a quadrilateral is unique; in particular, it is isomorphic to its dual quadrilateral $\begin{array}{ccc}
\bar{P} & \subset& M_1 \\
\cup & & \cup \\
M &\subset &\bar{Q} 
\end{array}$, where $$N \subset M \subset M_1, P \subset M \subset \bar{P}, Q \subset M \subset \bar{Q} $$ are each the basic construction. There is also present another intermediate subfactor $N \subset R \subset M$ with index $[M:R]=2 $.

In this quadrilateral $P$ and $Q$ are inner conjugate, and in fact a specific unitary $u \in N' \cap M_1$ which conjugates $\bar{P} $ onto $\bar{Q} $ is given by $u=2e_R-1 $, where $e_R$ is the biprojection in $M_1$ associated to $R$ \cite[Corollary 7.3.4]{G.J}.

We shall require the following results from \cite{G.J}, which hold for any two intermediate subfactors $P$ and $Q$ of a finite index irreducible inclusion $N\subseteq M$:
\begin{thm}\label{Landau}(Landau) $\displaystyle e_P\circ e_Q=\frac{\tr(e_Pe_Q)}{\delta}e_{PQ}$
\end{thm}
where $PQ$ is the (necessarily) strongly-closed subspace of $M$ generated by sums of products of elements in $P$ followed by elements of $Q$ and $\delta $ is $[M:N]^{\frac{1}{2}}$.  
\begin{prop}\label{traces} $\tr(e_{PQ})\tr(e_Pe_Q)=\tr(e_P)\tr(e_Q).$
\end{prop}
Using these facts, we can then prove:
\begin{lem}
 Let $\begin{array}{ccc}
P & \subset& M \\
\cup & & \cup \\
N &\subset &Q 
\end{array}$ be the unique irreducible quadrilateral of hyperfinite II$_1$ factors such that the index of each elementary subfactor is $2+\sqrt{2} $. There exists a unitary $v \in M $ such that $vPv^*=Q $ and $E_{P}(v)=0$.
\end{lem}
\begin{proof}
 By self-duality of the quadrilateral, it suffices to show that $E_{\bar{P}}(u)=0 $ for the $u$ described above. Let $e_{\bar{P}} \in M_2$ be the biprojection associated to $\bar{P} $. We then have $E_{\bar{P}}(u)e_{\bar{P}}=e_{\bar{P}}ue_{\bar{P}} $. It suffices to show that this quantity is zero, which we will do using the pictorial calculus of planar algebras.

Recall that an element $x$ in the relative commutant $N' \cap M_1 $ is represented by the  ``$2$-box'' \vpic{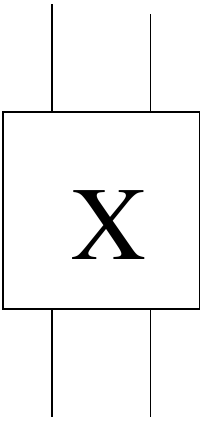}{.3in}, and more generally, an element of $N' \cap M_k $ by the ``$k$-box'' \vpic{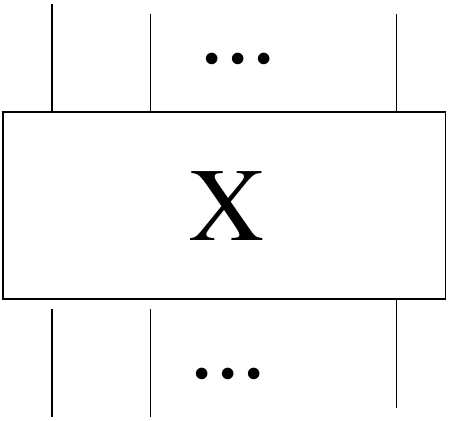}{.65in}. For details on planar algebras, see \cite{Jones.PA}. (Note that we follow the convention of \cite{G.J}, omitting the ``outer boundary''.) If \vpic{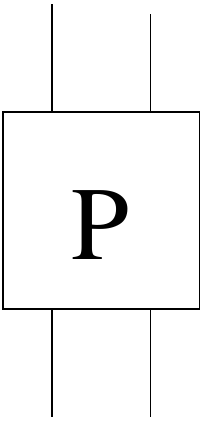}{.3in} is the $2$-box representing $e_P$, then by \cite[Lemma 3.2.6]{G.J} $e_{\bar{P}}=\displaystyle \frac{\delta}{\tr(e_P)}\phi(e_P)=\displaystyle \frac{\delta}{\tr(e_P)} $\vpic{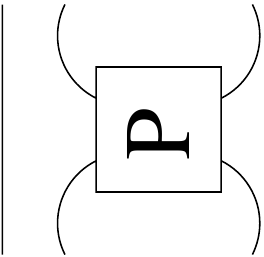}{.5in}$\in M'\cap M_2$ , where the modulus of the planar algebra $\delta $ is given by $[M:N]^{\frac{1}{2}}=2+\sqrt{2}=\frac{\tr(e_P)}{\delta} $ and $\phi$ is the linear isomorphism $x\to\delta^3\ce{M'}{xe_1e_2}$. Using Bisch's exchange relation for biprojections \cite{Bisch.Note} 
\begin{center}

\vpic{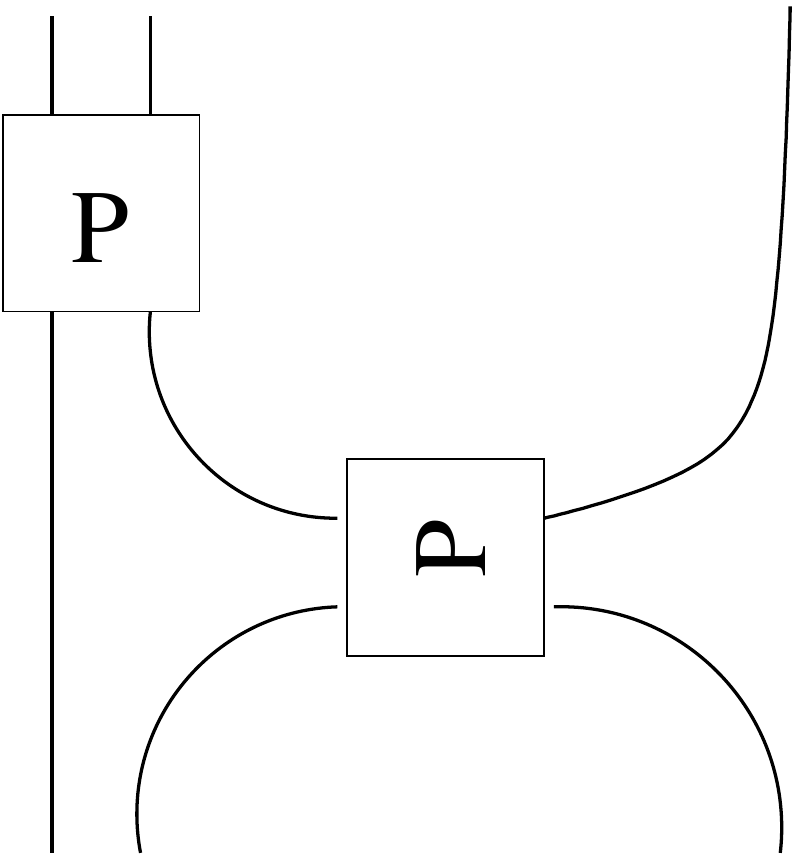} {1 in}  $=$ \vpic{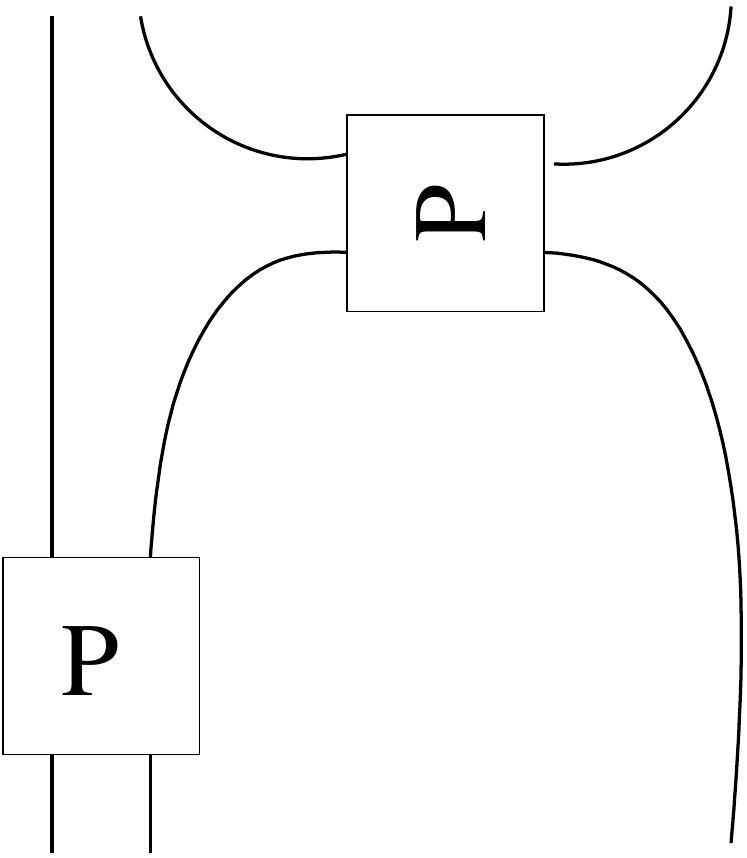} {1 in}

\end{center}

we compute:

\begin{center}
\vpic{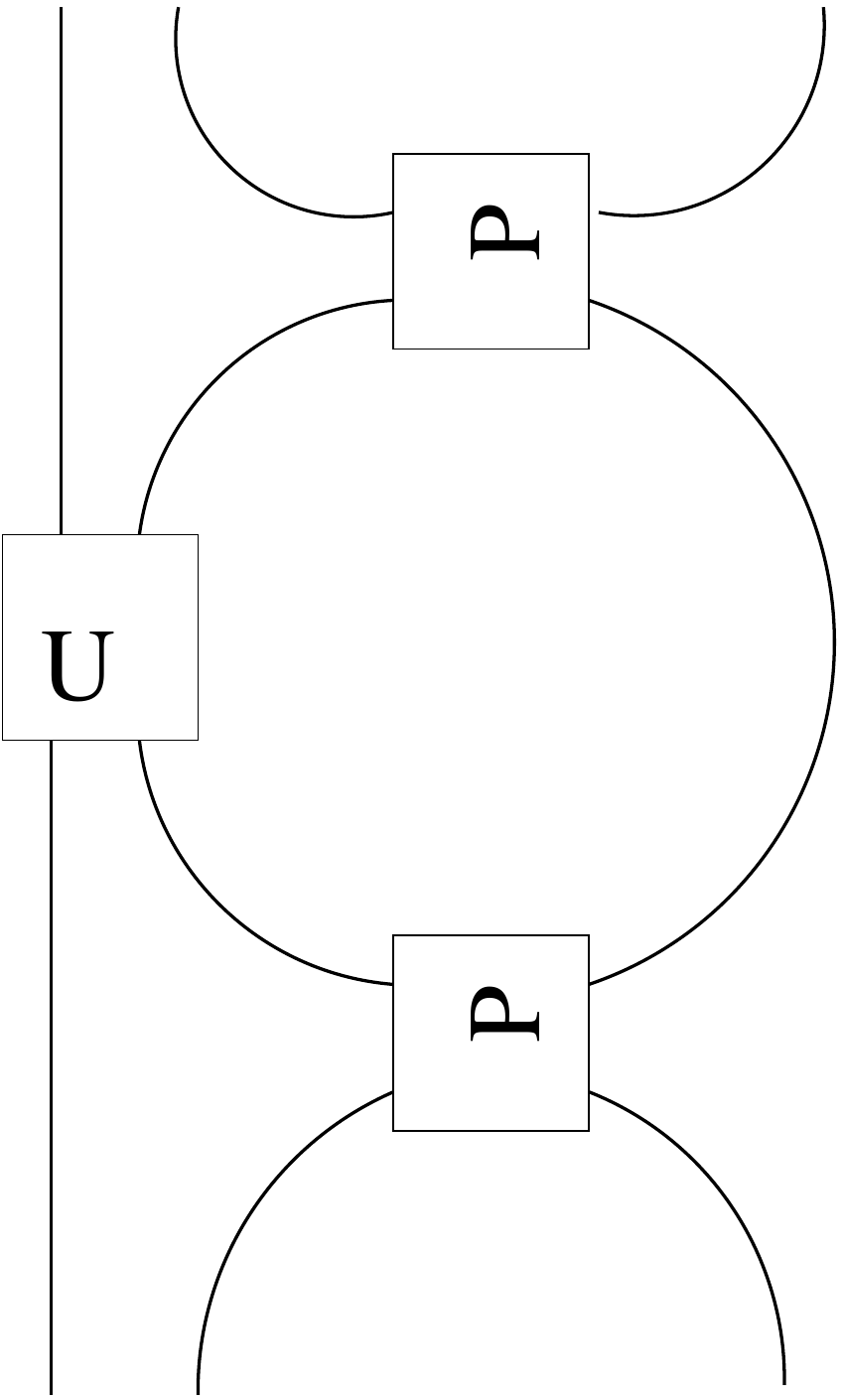}{1in}=\vpic{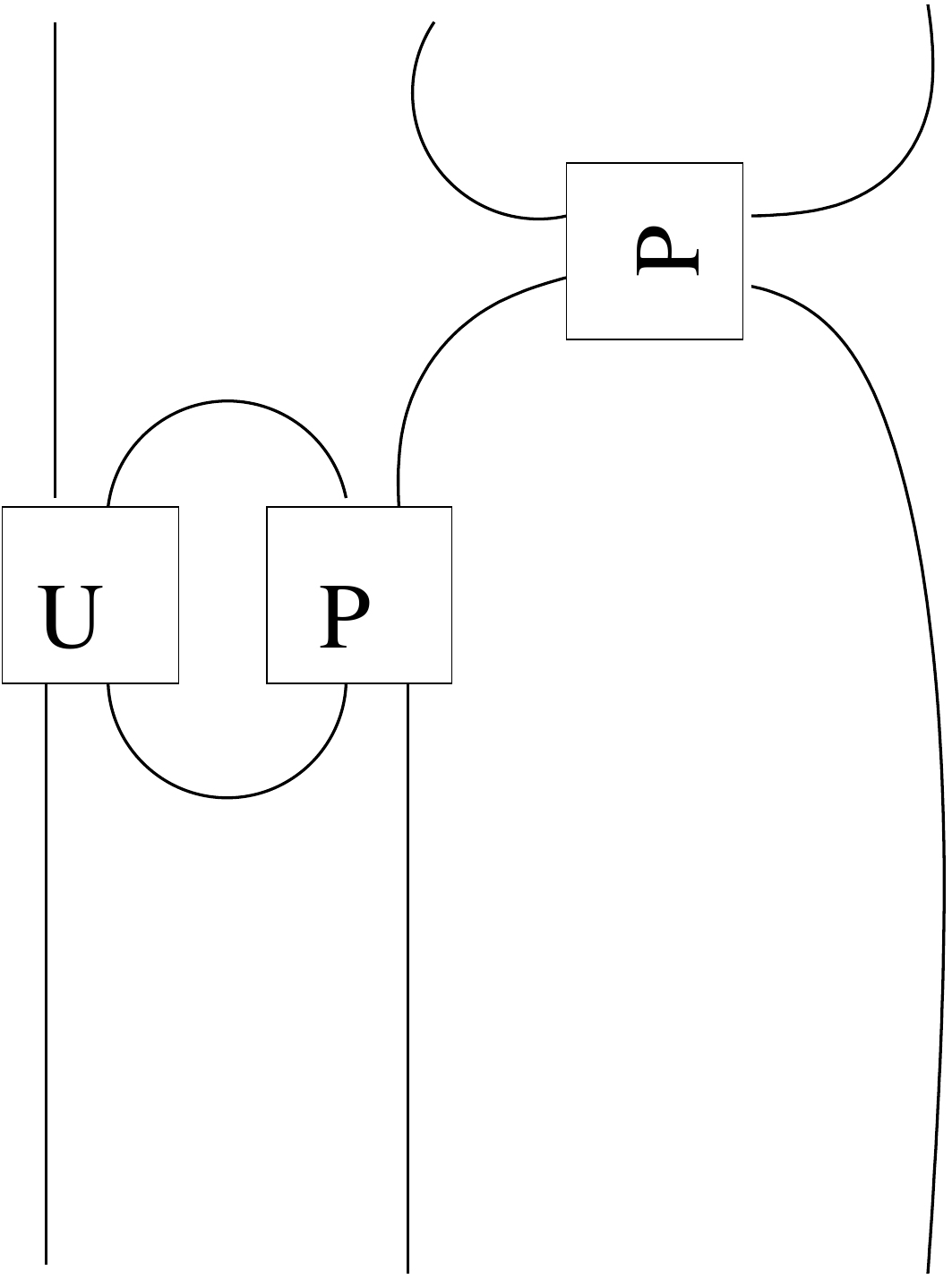}{1.3in}
\end{center}

$=e_{\bar{P}}\cdot (u \circ e_P)$, where $\circ $ is the comultiplication in the planar algebra. By Theorem \ref{Landau} and Proposition \ref{traces},
\begin{align*}
u \circ e_P=2e_R \circ e_P - 1 \circ e_P&=\frac 1 {\delta} \left(2\tr(e_Pe_R)e_{RP}-\tr(e_P)I\right)\\
&=\frac 1 {\delta}\left(2\left(\frac{\tr(e_R)\tr(e_P)}{\tr(e_{RP})}\right)e_{RP}-\tr(e_P)I\right).
\end{align*}
Since $R$ and $P$ cocommute (\cite[Lemma 7.3.1]{G.J}), by \cite[Theorem 3.2.1]{G.J} we have $RP=MP=M$ and so $e_{RP}=I$. Finally, $\tr(e_R)=\frac{\delta^2}2$, so that
\[
u \circ e_P=\frac 1 {\delta}\left(2\frac{\delta^2\tr(e_P)}{2\delta^2}-\tr(e_P)\right)I=0,
\]
therefore $E_{\bar{P}}(u)=0$ as well.
\end{proof}

\begin{lem}
 Let $x \in M $ with $||x||_2=1$. Then $||\ce{P}{x}-\ce{Q}{x}||_2  \leq \sqrt{2\sqrt{2}-2} $.
\end{lem}

\begin{proof}
 Consider the decomposition of $M$ into $N-N$ bimodules, which can be computed from the principal graph of $N \subset M $. We have $M \cong N \oplus 2V_1 \oplus 2V_2 \oplus V_3 $, where $N, V_1, V_2, V_3 $ are distinct irreducible $N-N$ bimodules (reference). Moreover, we have  $P \cong Q \cong N \oplus V_1 $. Without loss of generality, we may assume that $x$ lies in the $V_1 \oplus V_1$ component of $M$, since any component of $x$ contained in $N$ would be preserved under both $E_P$ and $E_Q$ and any component of $x$ contained in $V_2 \oplus V_2 \oplus V_3$ would vanish under both $E_P$ and $E_Q$. 

On the Hilbert space $L^2(V_1) \oplus L^2(V_1) $, the operators $e_P$ and $e_Q$ can be represented (after a suitable unitary transformation) in block form as the $2\times2$ matrices $A=\left( \begin{array}{cc}
1 & 0  \\
0 & 0  \end{array} \right)$ and $B=\left( \begin{array}{cc}
\lambda & \sqrt{\lambda(1-\lambda)}  \\
\sqrt{\lambda(1-\lambda)} & 1-\lambda  \end{array} \right)$, where the parameter $\lambda$ is the square of the cosine of the angle between these two projections. This angle is the same as the angle between the subfactors $P$ and $Q$, computed in \cite{G.J} as $Ang(P,Q)=\cos^{-1}(\sqrt{2}-1) $. So $\lambda=3-2\sqrt{2} $, and the norm of $e_P-e_Q $ on $L^2(V_1) \oplus L^2(V_1) $ is given by the positive eigenvalue of the matrix $A-B$, which is $\sqrt{1-\lambda}=\sqrt{2\sqrt{2}-2} $. 

\end{proof}
\begin{cor}
 $||\mathbb{E}_P-\mathbb{E}_{uPu^*}||_{\infty,2} \leq  \sqrt{2\sqrt{2}-2}$.
\end{cor}

\begin{proof}
If $x$ in $M$ satisfies $||x||_{\infty}\leq 1 $, then $||x||_{2} \leq ||x||_{\infty} \leq 1 $ 
so $||\mathbb{E}_P(x)-\mathbb{E}_Q(x)||_2  \leq  \sqrt{2\sqrt{2}-2} $.
\end{proof}

\begin{thm}\label{final}
 Let $A \subset M $ be the unique (up to conjugacy) subfactor of the hyperfinite II$_1$ factor with index $2+ \sqrt{2}$. Then $A$ is singular in $M$ but is not strongly singular.
\end{thm}

\begin{proof}
 Since $[M:A]=2+\sqrt{2}$ is between $3$ and $4$, $A \subseteq M$ is singular. By the previous two lemmas, the inequality \ref{ss} does not hold for $\alpha=1 $ and the $u$ described above, so the subfactor is not strongly singular. 
\end{proof}
Note that by combining Theorem \ref{final} with Theorem \ref{relativecomm}, we obtain that the optimal value for $\alpha$ in equation \eqref{ss} is between $\sqrt{\sqrt{2}(\sqrt{2}-1)}$ and $\sqrt{2(\sqrt{2}-1)}$.

\section{The WAHP and finite index singular subfactors}\label{nowahp}

A subalgebra $B$ is said to have the weak asymptotic homomorphism property (WAHP) if for every $\varepsilon >0$ and for all $x_1,\dots, x_n$, $y_1,\dots, y_m$ in $M$, there exists a unitary $u$ in $B$ with
\begin{equation}\label{wahp}
\n \ce{B}{x_iuy_j}-\ce{B}{x_i}u\ce{B}{y_j}\n_2<\varepsilon
\end{equation}
for every $1\leq i\leq n$, $1\leq j\leq m$. 

As a consequence of Popa's Intertwining Theorem \cite{Popa.StrongrigidityI}, a subalgebra $B$ of a \IIi factor $M$ will have the WAHP if and only if there are no nonzero finite projections in $B'\cap \langle M,e_B\rangle$ subordinate to $e_B^{\perp}$. We include a simple proof showing that no finite index subfactor may have the WAHP. 

Before beginning the proof, recall from \cite{Jones.SubfactorsBook} or \cite{Popa.Entropy} that a Pimsner-Popa basis for $N\subseteq M$ a finite index inclusion of subfactors is a collection of elements $\lambda_1, \dots ,\lambda _k$ in $M$ with $k$ any integer greater than or equal to $[M:N]$ such that every $x\in M$ may be represented as $\sum^k _{j=1} \lambda_j \ce{N}{\lambda^* _jx}$ and $\sum^k _{j=1} \lambda _j e_N\lambda _j^*=1$.

\begin{thm}\label{nwahp}
If $N\subseteq M$ is a {\rm{\IIi}}factor with $1<[M:N]<\infty$, then $N$ does not have the WAHP in $M$.
\end{thm}
\begin{proof} It will be advantageous to use a Pimsner-Popa basis obtained by first choosing $k$ to be the least integer greater than or equal to $[M:N]$. We then select a collection of orthogonal projections $\{p_j\}_{j=1}^k$ in $\langle M,e_N\rangle$ with $p_1=e_N$, $\displaystyle \sum_{j=1} ^k p_j=1$ and ${\textrm Tr}(p_i)\leq 1$ with equality except possibly for $j=k$. 

Let $v_1,v_2,\dots ,v_k$ be partial isometries in $\langle M,e_N\rangle$ such that $e_N=v_1$, $v_jv^* _j=p_j$ and $v_j ^*v_j\leq e_N$. The desired basis is given by the unique elements $\lambda_j\in M$ with the property that 
\[
\lambda_je_N=v_je_N.
\]
Observe that $\lambda_1=1$. Since for $i\ne j$, 
\begin{align*}
\ce{N}{\lambda_i^* \lambda_j}e_N&=e_N\lambda_i^*\lambda_je_N=e_Nv_i^*v_je_N\\
&=e_Nv_i^*p_ip_jv_je_N=0,
\end{align*}
we have that $\ce{N}{\lambda_i^*\lambda_j}=0$ for $i\ne j$. In particular, $\ce{N}{\lambda_j}=0$ for all $1<j\leq k$. It is worth noting that this is the original construction in \cite{Popa.Entropy}.

Now suppose $1<[M:N]<\infty$ and $\lambda_1,\dots ,\lambda_k$ are chosen as indicated. We will show that the WAHP fails for the sets $\{ x_i=\lambda_i\}$ and $\{ y_j=\lambda_j^*\}$, $1<i,j\leq k$. Let $u$ be any unitary in $N$. Then since 
\[
\tau (\ce{N}{x})=\tau (x)=\tr(e_Nx)
\]
for all $x$ in $M$,
\begin{align*}
\sum_{i,j=2} ^k \n \ce{N}{\lambda_i ^*u\lambda_j}\n^2 _2&=\sum_{i,j=2} ^k \tau(\ce{N}{\lambda_j^* u^*\lambda_i}\ce{N}{\lambda_i^* u\lambda_j})=\sum_{i,j=2} ^k \tau(\lambda_j^* u^*\lambda_i\ce{N}{\lambda_i^* u\lambda_j})\\
&=\sum_{i,j=2} ^k \tr (e_N\lambda_j^* u^*\lambda_i\ce{N}{\lambda_i^* u\lambda_j})=\sum_{i,j=2} ^k \tr (e_N\lambda_j^* u^*\lambda_ie_N{\lambda_i^* u\lambda_j}e_N).
\end{align*}
Using this equality, the fact that $u$ commutes with $e_N$, and $\sum_{j=1}^k\lambda_je_N\lambda_j^*=1$, we get
\begin{align*}
\sum_{i,j=2} ^k \n \ce{N}{\lambda_i ^*u\lambda_j}\n^2 _2&=\sum_{i,j=2} ^k \tr (e_N\lambda_j^* u^*\lambda_ie_N{\lambda_i^* u\lambda_j}e_N)=\sum_{i,j=2} ^k \tr(u^*\lambda_ie_N\lambda_i^*u\lambda_je_N\lambda_j^*)\\
&=\tr (u^*(1-e_N)u(1-e_N))=\tr ((1-e_N)u^*u)=[M:N]-1>0.
\end{align*}
This implies that for any given unitary $u$ in $N$, there are indices $1<i,j\leq k$ with 
\[
 \n \ce{N}{\lambda^* _iu\lambda_j}\n_2\geq \frac {\sqrt{[M:N]-1}}{k-1},
 \]
and so the WAHP fails to hold. 
\end{proof}

Combining the previous theorem with the discussion at the beginning of this section, we immediately get
\begin{cor} 
There exist singular subfactors that do not have the WAHP.
\end{cor}

\section{Acknowledgements}
The second named author would like to thank his advisor, Professor Roger Smith, for guidance and support and Professors Stuart White and Allan Sinclair for helpful discussions. The authors would also like to thank Professor Vaughan Jones for suggesting the examination of the $A_7$ quadrilateral and Professor Dietmar Bisch for helpful conversations. Portions of this work are taken from the second named author's dissertation, completed under the supervision of Professor Smith.

\bibliographystyle{plain}

\begin{bibdiv}
\begin{biblist}

\bib{Dietmar.Index6}{article}{
      author={Bisch, Dietmar},
      author={Nicoara, Remus},
      author={Popa, Sorin},
       title={Continuous families of hyperfinite subfactors with the same
  standard invariant},
        date={2007},
        ISSN={0129-167X},
     journal={Internat. J. Math.},
      volume={18},
      number={3},
       pages={255\ndash 267},
      review={\MR{MR2314611 (2008k:46188)}},
}

\bib{Bisch.Note}{article}{
      author={Bisch, Dietmar},
       title={A note on intermediate subfactors},
        date={1994},
        ISSN={0030-8730},
     journal={Pacific J. Math.},
      volume={163},
      number={2},
       pages={201\ndash 216},
      review={\MR{MR1262294 (95c:46105)}},
}

\bib{Dixmier.Masa}{article}{
      author={Dixmier, J.},
       title={Sous-anneaux ab\'eliens maximaux dans les facteurs de type fini},
        date={1954},
        ISSN={0003-486X},
     journal={Ann. of Math. (2)},
      volume={59},
       pages={279\ndash 286},
      review={\MR{MR0059486 (15,539b)}},
}

\bib{Fuglede.MaxInjective}{article}{
      author={Fuglede, Bent},
      author={Kadison, Richard~V.},
       title={On a conjecture of {M}urray and von {N}eumann},
        date={1951},
     journal={Proc. Nat. Acad. Sci. U. S. A.},
      volume={37},
       pages={420\ndash 425},
      review={\MR{MR0043390 (13,255a)}},
}

\bib{Goldman.Index2}{article}{
      author={Goldman, M.},
       title={On subfactors of factors of type {${\rm II}\sb 1$}},
        date={1960},
     journal={Mich. Math. J.},
      volume={7},
       pages={167\ndash 172},
}

\bib{G.J}{article}{
      author={Grossman, Pinhas},
      author={Jones, Vaughan F.~R.},
       title={Intermediate subfactors with no extra structure},
        date={2007},
        ISSN={0894-0347},
     journal={J. Amer. Math. Soc.},
      volume={20},
      number={1},
       pages={219\ndash 265 (electronic)},
      review={\MR{MR2257402 (2007h:46077)}},
}

\bib{Jones.SubfactorsBook}{book}{
      author={Jones, V.},
      author={Sunder, V.~S.},
       title={Introduction to subfactors},
      series={London Mathematical Society Lecture Note Series},
   publisher={Cambridge University Press},
     address={Cambridge},
        date={1997},
      volume={234},
        ISBN={0-521-58420-5},
      review={\MR{MR1473221 (98h:46067)}},
}

\bib{Jones.Index}{article}{
      author={Jones, V. F.~R.},
       title={Index for subfactors},
        date={1983},
        ISSN={0020-9910},
     journal={Invent. Math.},
      volume={72},
      number={1},
       pages={1\ndash 25},
      review={\MR{84d:46097}},
}

\bib{Jones.PA}{unpublished}{
      author={Jones, V.F.R.},
       title={Planar algebras, i},
        date={1999},
        note={Preprint, arXiv:math/9909027v1},
}

\bib{Popa.Entropy}{article}{
      author={Pimsner, Mihai},
      author={Popa, Sorin},
       title={Entropy and index for subfactors},
        date={1986},
        ISSN={0012-9593},
     journal={Ann. Sci. \'Ecole Norm. Sup. (4)},
      volume={19},
      number={1},
       pages={57\ndash 106},
      review={\MR{MR860811 (87m:46120)}},
}

\bib{Sinclair.PertSubalg}{article}{
      author={Popa, Sorin},
      author={Sinclair, Allan~M.},
      author={Smith, Roger~R.},
       title={Perturbations of subalgebras of type {II{$\sb 1$}} factors},
        date={2004},
        ISSN={0022-1236},
     journal={J. Funct. Anal.},
      volume={213},
      number={2},
       pages={346\ndash 379},
      review={\MR{MR2078630}},
}

\bib{Popa.Kadison}{article}{
      author={Popa, Sorin},
       title={On a problem of {R}. {V}. {K}adison on maximal abelian {$\ast
  $}-subalgebras in factors},
        date={1981/82},
        ISSN={0020-9910},
     journal={Invent. Math.},
      volume={65},
      number={2},
       pages={269\ndash 281},
      review={\MR{MR641131 (83g:46056)}},
}

\bib{Popa.MaxInjective}{article}{
      author={Popa, Sorin},
       title={Maximal injective subalgebras in factors associated with free
  groups},
        date={1983},
        ISSN={0001-8708},
     journal={Adv. in Math.},
      volume={50},
      number={1},
       pages={27\ndash 48},
      review={\MR{MR720738 (85h:46084)}},
}

\bib{Popa.SingularMasas}{article}{
      author={Popa, Sorin},
       title={Singular maximal abelian {$\ast $}-subalgebras in continuous von
  {N}eumann algebras},
        date={1983},
        ISSN={0022-1236},
     journal={J. Funct. Anal.},
      volume={50},
      number={2},
       pages={151\ndash 166},
      review={\MR{MR693226 (84e:46065)}},
}

\bib{Popa.CBMSNotes}{book}{
      author={Popa, Sorin},
       title={Classification of subfactors and their endomorphisms},
      series={CBMS Regional Conference Series in Mathematics},
   publisher={Published for the Conference Board of the Mathematical Sciences,
  Washington, DC},
        date={1995},
      volume={86},
        ISBN={0-8218-0321-2},
      review={\MR{MR1339767 (96d:46085)}},
}

\bib{Popa.StrongrigidityI}{article}{
      author={Popa, Sorin},
       title={Strong rigidity of {$\rm II\sb 1$} factors arising from malleable
  actions of {$w$}-rigid groups. {I}},
        date={2006},
        ISSN={0020-9910},
     journal={Invent. Math.},
      volume={165},
      number={2},
       pages={369\ndash 408},
      review={\MR{MR2231961 (2007f:46058)}},
}

\bib{Sinclair.StrongSing2}{article}{
      author={Robertson, Guyan},
      author={Sinclair, Allan~M.},
      author={Smith, Roger~R.},
       title={Strong singularity for subalgebras of finite factors},
        date={2003},
        ISSN={0129-167X},
     journal={Internat. J. Math.},
      volume={14},
      number={3},
       pages={235\ndash 258},
      review={\MR{2004c:22007}},
}

\bib{Sinclair.strongsing}{article}{
      author={Sinclair, A.~M.},
      author={Smith, R.~R.},
       title={Strongly singular masas in type {$\rm II\sb 1$} factors},
        date={2002},
        ISSN={1016-443X},
     journal={Geom. Funct. Anal.},
      volume={12},
      number={1},
       pages={199\ndash 216},
      review={\MR{2003i:46066}},
}

\bib{singular.al}{article}{
      author={Sinclair, Allan~M.},
      author={Smith, Roger~R.},
      author={White, Stuart~A.},
      author={Wiggins, Alan},
       title={Strong singularity of singular masas in {${\rm II}\sb 1$}
  factors},
        date={2007},
        ISSN={0019-2082},
     journal={Illinois J. Math.},
      volume={51},
      number={4},
       pages={1077\ndash 1084},
      review={\MR{MR2417416}},
}

\bib{Vaes.nontrivialsubfactor}{unpublished}{
      author={Vaes, Stefaan},
       title={Factors of type {II{$\sb 1$}} without non-trivial finite index
  subfactors},
        date={2006},
        note={Preprint, arXiv:math.OA/0610231},
}

\end{biblist}
\end{bibdiv}

  
\section*{Author Addresses}

\begin{tabular*}{\textwidth}{l@{\hspace*{2cm}}l}
 Pinhas Grossman&Alan Wiggins\\
   Department of Mathematics&Department of Mathematics\\
   1326 Stevenson Center&1326 Stevenson Center\\
   Vanderbilt University&Vanderbilt University\\ 
   Nashville, TN, 37209&Nashville, TN, 37209\\
   USA&USA\\
\texttt{pinhas.grossman@vanderbilt.edu}&\texttt{alan.d.wiggins@vanderbilt.edu}
{}\\\\
\end{tabular*}

\end{document}